\documentstyle[12pt]{article}
\newcommand{\Rset}[0]{\mbox{\rm I\kern-.200em R}}
\begin{document}
\centerline{\bf A RIEMANN SUM UPPER BOUND IN THE}
\centerline{\bf RIEMANN-LEBESQUE THEOREM}
\mbox{}\\
\mbox{}\hskip1.5in MAURICE H.P.M. VAN PUTTEN
	\footnote{Department of 
	Mathematics, 
	MIT, Cambridge, MA 02139.}
\mbox{}\\
\mbox{}\\
{\small{\bf Abstract.} The Riemann-Lebesque Theorem is commonly proved in a few
strokes using the theory of Lebesque integration.
Here, 
the upper bound 
$2\pi|c_k(f)|\le S_k(f)-s_k(f)$ for the Fourier
coefficients $c_k$ is proved
in terms of majoring and minoring Riemann
sums $S_k(f)$ and $s_f(k)$, respectively, for Riemann
integrable functions $f(x)$. This proof has been used
in a course on methods of applied
mathematics.}\\
\mbox{}
\mbox{}
\mbox{}\\
One of the central theorems in the theory of Fourier
series is the 
Riemann-Lebesque theorem: if 
\begin{eqnarray}
c_k=\frac{1}{2\pi}\int_0^{2\pi}
f(x)e^{ikx}\mbox{d}x
\end{eqnarray}
are the Fourier coefficients of an integrable
function $f(x)$ on $[0,2\pi]$, then
\begin{eqnarray}
c_k\rightarrow0
\end{eqnarray}
as $k$ becomes large. In courses on methods of applied mathematics,
one is usually restricted to present the theorem in the context
of Riemann integrability. While intuition makes the result plausible,
traditional proofs tend to be 
somewhat involved.

In the present proof, we exploit the length scale $\lambda=\frac{2\pi}{k}$
introduced on
$[0,2\pi]$ following a choice of $k\epsilon {\bf N}$.
Put $x_n=n\lambda$, $n=0,\cdots,k$. Then
\begin{eqnarray}
2\pi c_k=\int_0^{2\pi}f(x)e^{ikx}
\mbox{d}x=\Sigma_{n=0}^{k-1}\int_{x_n}^{x_{n+1}}
f(x)e^{ikx}\mbox{d}x.
\end{eqnarray}
Focusing on a particular subinterval $[x_n,x_{n+1}]$, set
$x=x_n+\frac{\lambda}{2\pi}y$, $y\epsilon[0,2\pi]$. Thus,
$e^{ikx}=e^{ik(x_n+\frac{\lambda}{2\pi}y)}=e^{iy}$, and so
\begin{eqnarray}
\int_{x_n}^{x_{n+1}}f(x)e^{ikx}\mbox{d}x=
\frac{\lambda}{2\pi}\int_0^{2\pi}f(x_n+\frac{\lambda}{2\pi}y)
e^{iy}\mbox{d}y
\end{eqnarray}
Following `adding and subtraction,' we obtain
\begin{eqnarray}
\begin{array}{rll}
\int_0^{2\pi}f(x_n+\frac{\lambda}{2\pi}y)\mbox{d}y&=&
\int_0^{2\pi}[f(x_n+\frac{\lambda}{2\pi}y)-f(x_n)]
e^{iy}\mbox{d}y
+\int_0^{2\pi}f(x_n)e^{iy}\mbox{d}y\\
&=&\int_0^{2\pi}[f(x_n+\frac{\lambda}{2\pi}y)-f(x_n)]
e^{iy}\mbox{d}y,
\end{array}
\end{eqnarray}
since integration of a constant against $e^{iy}$ over
a full period is identically zero.
Upon taking absolute values, it follows that
\begin{eqnarray}
2\pi |c_k|\le \Sigma_{n=0}^{k-1}\frac{\lambda}{2\pi}
\int_0^{2\pi}|[f(x_n+\frac{\lambda}{2\pi}y)
-f(x_n)]|\mbox{d}y.
\end{eqnarray}
To finalize, introduce the majoring
and minoring Riemann sums
\begin{eqnarray}
S_k(f)=\lambda\Sigma_{n=0}^{k-1} f_{kn}^S,
\mbox{  }
s_k(f)=\lambda\Sigma_{n=0}^{k-1} f_{kn}^I,
\end{eqnarray}
where $f_{kn}^S=\mbox{sup}_{[x_n,x_{n+1}]}f(x)$ and
$f_{kn}^I=\mbox{inf}_{[x_n,x_{n+1}]}f(x)$.
Then 
\begin{eqnarray}
|[f(x_n)+\frac{\lambda}{2\pi}y)-
f(x_n)]|\le f^S_{kn}-f^I_{kn}
\end{eqnarray}
on $[x_n,x_{n+1}]$, so that
\begin{eqnarray}
2\pi |c_k|\le \Sigma_{n=0}^{k-1}\lambda(f^S_{kn}-f^I_{kn})=S_k(f)-s_k(f).
\end{eqnarray}
Clearly, an assumption of 
Riemann integrability assures
$c_k\rightarrow0$ as $k\rightarrow\infty$.
This completes the proof.

Riemann integrability
enables the preceeding arguments 
to be readily generalized to 
the case of an infinite interval,
and, similarly, 
to the case of non-integer
$k$. However, in our experience
a preferred presentation in
an introductory course on applied mathematics courses
($e.g.$ [1]) is 
in the context of an
interval $[0,2\pi]$.\\
\mbox{}\\
\centerline{\small REFERENCES}

 [1] G. Strang,{\em Introduction 
 to Applied Mathematics},
 Wellesley-Cambridge Press,
 Cambridge, MA, 1996.
\end{document}